\documentclass[12pt]{amsart}
\usepackage{somedefs,tracefnt,latexsym,,rawfonts,epsfig,amsfonts,amssymb,latexsym,enumerate,amscd}
\def\cal{\mathcal}
\def\Bbb{\mathbb}

\def\G{\Gamma}
\def\r{\rangle}
\def\l{\langle}

\def\ul{\underline}
\def\F{\cal {\ul{FICWF}}}
\def\uF{\cal {FICWF}}
\def\flf{FICwF^{_{\underline L}{\cal H}^?_*}_{\cal {FIN}}}
\def\flv{FICwF^{_{\underline L}{\cal H}^?_*}_{\cal {VC}}}

\def\ficf{FIC^{{\cal H}^?_*}_{\cal {FIN}}}
\def\ficwf{FICwF^{{\cal H}^?_*}_{\cal {FIN}}}
\def\ficlu{FIC^{_{L}{\cal H}^?_*}_{\cal {C}}}
\def\ficl{FIC^{_{\underline L}{\cal H}^?_*}_{\cal {C}}}
\def\ficwl{FICwF^{_{\underline L}{\cal H}^?_*}_{\cal {C}}}
\def\ficwlu{FICwF^{_{L}{\cal H}^?_*}_{\cal {C}}}

\def\flvu{FICwF^{_{L}{\cal H}^?_*}_{\cal {VC}}}

\def\flvuw{FIC^{_{L}{\cal H}^?_*}_{\cal {VC}}}
\def\flvw{FIC^{_{\ul L}{\cal H}^?_*}_{\cal {VC}}}
\def\flfw{FIC^{_{\ul L}{\cal H}^?_*}_{\cal {FIN}}}
\def\pvc{{\cal P}^{_L{\cal H}^?_*}_{\cal {VC}}}
\def\pvcu{{\cal P}^{_{\ul L}{\cal H}^?_*}_{\cal {VC}}}
\def\pfu{{\cal P}^{_{\ul L}{\cal H}^?_*}_{\cal {FIN}}}

\def\wtfinu{_w{\cal T}^{_{\ul L}{\cal H}^?_*}_{\cal {FIN}}}

\def\wttvu{_{wt}{\cal T}^{_{\ul L}{\cal H}^?_*}_{\cal {VC}}}
\def\wttv{_{wt}{\cal T}^{_{L}{\cal H}^?_*}_{\cal {VC}}}
\def\tfinv{{\cal T}^{_{L}{\cal H}^?_*}_{{\cal {FIN}}, {\cal {VC}}}}
\def\tfinvu{{\cal T}^{_{\ul L}{\cal H}^?_*}_{{\cal {FIN}}, {\cal {VC}}}}

\newtheorem{prop}{Proposition}[section]
\newtheorem{thm}{Theorem}[section]

\newtheorem*{mainl}{Main Lemma}
\newtheorem{lemma}{Lemma}[section]
\newtheorem{cor}{Corollary}[section]

\newtheorem{defn}{Definition}[section]

\newtheorem{rem}{Remark}[section]
\numberwithin{equation}{section}
\begin{document}
\date{July 03, 2008.}
\title[The fibered isomorphism conjecture in $L$-theory]{The 
fibered isomorphism conjecture in $L$-theory}
\author[S.K. Roushon]{S. K. Roushon}
\address{School of Mathematics\\
Tata Institute\\
Homi Bhabha Road\\
Mumbai 400005, India}
\email{roushon@math.tifr.res.in} 
\urladdr{http://www.math.tifr.res.in/\~\ roushon/}
\begin{abstract} This is the first of three articles on the 
Fibered Isomorphism Conjecture of Farrell and Jones for $L$-theory. 
We apply the general techniques developed in \cite{R1} and 
\cite{R2} to the $L$-theory case of the 
conjecture and prove several 
results.

Here we prove the conjecture, 
after inverting $2$, for poly-free groups. In particular, 
it follows for braid groups. 
We also prove the conjecture for some classes of groups 
without inverting $2$.
In fact we consider a general class of groups satisfying 
certain conditions  
which includes the above groups and some other important 
classes of groups. We check that 
the properties we  
defined in \cite{R1} are satisfied in several instances 
of the conjecture. 

\end{abstract} 



\maketitle

%\tableofcontents

\section{Introduction} This is the first of three articles 
where we study the Fibered Isomorphism Conjecture of Farrell 
and Jones ([\cite{FJ}, 1.7]) for 
$L^{\l -\infty \r}$-theory. The Isomorphism Conjecture 
was stated in [\cite{FJ}, 1.6] for the 
pseudoisotopy theory, $K$-theory and the $L^{\l -\infty \r}$-theory. 
The conjecture says that these theories can be 
computed for a group if we can compute them for all 
its virtually cyclic subgroups. It is known that the 
conjecture is not true for the other $L$-theory functors 
$L^h$ and $L^s$ (see \cite{FJL}). The Fibered Isomorphism  
Conjecture is stronger and is appropriate 
for induction arguments. This property is crucial for the 
method we use here.
We developed some general techniques in \cite{R1} 
and \cite{R2} which applies to all the three cases 
of the conjecture and proved the Fibered Isomorphism 
Conjecture for 
the pseudoisotopy case for several classes of groups. 

Here we see that the techniques can be used effectively 
to prove the Fibered Isomorphism Conjecture for 
$L^{\l -\infty \r}$-theory and for $\ul L^{\l -\infty \r}$-theory for some 
well-known classes of groups, where ${\ul L}^{\l -\infty \r}={
L}^{\l -\infty \r}\otimes {\Bbb Z}[\frac{1}{2}]$. The advantage 
of taking $\ul L^{\l -\infty \r}$-theory is that we can 
consider the `finite subgroups' instead of the `virtually cyclic 
subgroups'. This will also help to prove the conjecture for 
$\ul L^{\l -\infty \r}$-theory for a larger class of groups.

We also show that the two properties 
(${\cal T}_{\cal C}^{{\cal H}^?_*}$ and 
${\cal P}_{\cal C}^{{\cal H}^?_*}$) defined in \cite{R1} 
are satisfied in several instances of the conjecture 
for $L^{\l -\infty \r}$-theory. 

Throughout the article by `group' we mean `discrete countable group' 
unless otherwise mentioned.

\begin{defn} \label{begin} {\rm Let $\F$ ($\uF$) be the smallest
class of groups satisfying the conditions $1$ to 
$5$ (i to iv) below.

1. $\F$ contains the cocompact
discrete subgroups of linear Lie groups with finitely many
components.

2. (Subgroup) If $H<G\in \F$ then $H\in \F$

3. (Free product) If $G_1, G_2\in \F$ then $G_1*G_2\in \F$.

4. (Direct limit) If $\{G_i\}_{i\in I}$ is a directed sequence
of groups with
$G_i\in \F$. Then the limit $\lim_{i\in I}G_i\in \F$.

5. (Extension) For an exact sequence of groups $1\to K\to G\to N\to 1$, if
$K, N\in \F$ then $G\in \F$. 

i. $\uF$ contains the cocompact
discrete subgroups of Lie groups with finitely many
components. 

ii. 2, 3 and 4 as above after replacing $\F$ by $\uF$.

iii. (Direct product) If $G_1, G_2\in \uF$ then $G_1\times G_2\in \uF$.

iv. (Polycyclic extension) For an exact sequence of 
groups $1\to K\to G\to N\to 1$, if
either $K$ is virtually cyclic and  $N\in \uF$ or 
$N$ is finite and $K\in \uF$ then $G\in \uF$.}
\end{defn}

For two groups $A$ and $B$ the wreath product $A\wr B$ 
is the semidirect product $A^B\rtimes B$ 
with respect to the regular action of $B$ on $A^B$. Here $A^B$ 
denotes the 
direct sum of copies of $A$ indexed by $B$. 
Let ${\cal {VC}}$ and 
${\cal {FIN}}$ denote the class of virtually cyclic groups 
and the class
of finite groups respectively.

We prove the following Theorem. For notations see Section 2.

\begin{thm} \label{mth} Let $\Gamma\in \F$ and $\Delta\in \uF$. 
Let $G$ and $H$ be two groups with homomorphisms $\phi :G\to
\Gamma\wr F$ and $\psi :H\to \Delta\wr K$ where $F$ and 
$K$ are finite groups. 
Then the following assembly maps
are isomorphisms for all $n$. 
$${\cal
H}^G_n(E_{\phi^*{\cal
{FIN}}(\G\wr F)}(G), {\bf \ul L}^{\l -\infty \r})\to
{\cal H}^G_n(pt, {\bf \ul L}^{\l -\infty \r})\simeq
\ul L_n^{\l -\infty \r}({\Bbb
Z}G).$$
$${\cal
H}^H_n(E_{\psi^*{\cal
{VC}}(\Delta\wr K)}(H), {\bf L}^{\l -\infty \r})\to
{\cal H}^H_n(pt, {\bf L}^{\l -\infty \r})\simeq
 L_n^{\l -\infty \r}({\Bbb
Z}H).$$

In other words the Fibered Isomorphism Conjecture
of Farrell and Jones for the $\ul L^{\l -\infty \r}$-theory ( 
$L^{\l -\infty \r}$-theory)  
is true for the
group $\G\wr F$ ($\Delta\wr K$). Equivalently, 
the $\flf(\G)$  and the $\flvu(\Delta)$ are satisfied.\end{thm} 

The notations in the above statement are described in the next section.

\begin{thm} \label{mth1} Let $\cal C$ ($\cal D$) be the class of groups 
which satisfy the $\flf$ (the $\flvu$). 
Then $\cal C$ ($\cal D$) has the properties 2 to 5 
(ii to iv) after replacing $\F$ ($\uF$) by $\cal C$ 
($\cal D$) in Definition \ref{begin}.\end{thm}

Our next goal is to show that $\F$ and $\uF$ contains some 
well-known classes of groups.

\begin{thm} \label{hyper} $\F$ contains the following groups. 

1. Virtually cyclic groups.

2. Free groups and abelian groups.

3. Poly-free groups. A 
poly-free group $G$ admits a filtration by subgroups: 
$1<G_0<G_1<\cdots < G_n=G$ so that $G_i$ is normal in $G_{i+1}$ 
and $G_{i+1}/G_i$ is free. Here $n$ is called the index 
of $G$.

4. Strongly poly-free groups. See [\cite{AFR}, definition 1.1]. 

5. Full braid groups. 

6. Cocompact discrete subgroups of Lie 
groups with finitely many components. 

7. Groups whose some derived subgroup belong to $\F$.

$\uF$ contains the following groups. 

i. Virtually cyclic groups.

ii. Free groups and abelian groups.

iii. Groups appearing in 6 (by definition).
 
iv. Virtually polycyclic groups.\end{thm}   

\begin{rem}{\rm Here we should remark that the non-fibered version 
of the Isomorphism 
conjecture in $\ul L^{\l -\infty \r}$-theory for a class of groups 
including poly-free groups and one-relator groups  
was proved in [\cite{BL}, theorem 0.13].}\end{rem}

In the next section we recall the statement of the conjecture in 
the more general context in equivariant homology theory 
from \cite{BL}. Also we recall some notations and definitions 
from \cite{R1}. For the statement of the original Isomorphism 
Conjectures see [\cite{FJ}, 1.6, 1.7]. 

When $\Gamma$ and $\Delta$ are torsion free, $F=K=\{1\}$ and $\phi$ 
and $\psi$ are the
identity maps, Theorem \ref{mth} reduces to the isomorphism of the classical 
assembly map in surgery theory. Therefore we have the following 
corollary. See 1.6.1 in \cite{FJ} for details.

\begin{cor} \label{classic} Let $\G\in \F$ and $\Delta\in\uF$ 
and in 
addition assume that $\G$ and $\Delta$ are torsion free. Then the
following assembly maps are isomorphism for all $n$. $$H_n(B\Gamma, {\bf
\ul L}^{\l -\infty \r})\to \ul L^{\l -\infty \r}_n({\Bbb
Z}\Gamma).$$
$$H_n(B\Delta, {\bf
L}^{\l -\infty \r})\to L^{\l -\infty \r}_n({\Bbb
Z}\Delta).$$
\end{cor}

In other words the surgery groups $\ul L^{\l -\infty \r}_n({\Bbb Z}\Gamma)$ of
$\Gamma$ and the surgery groups $L^{\l -\infty \r}_n({\Bbb Z}\Delta)$ of
$\Delta$ form generalized homology theories.

Since the surgery groups with different 
decoration defer by $2$-torsions we also have the 
following. See [\cite{FL}, section 5, para 1]. 

\begin{lemma} \label{invert2} 
${\ul L}^{\l -\infty \r}_n({\Bbb Z}\Gamma)\simeq L^h_n({\Bbb
Z}\Gamma)\otimes {\Bbb Z}[\frac{1}{2}]\simeq 
L^s_n({\Bbb Z}\Gamma)\otimes {\Bbb Z}[\frac{1}{2}]$ 
for any group $\G$.\end{lemma}

Therefore, Lemma \ref{invert2} implies that the Theorem \ref{mth} 
is true for the functors $L^h\otimes {\Bbb Z}[\frac{1}{2}]$ 
and $L^s\otimes {\Bbb Z}[\frac{1}{2}]$ also. We have already 
mentioned in the Introduction that the 
Theorem \ref{mth} is not true for $L^h$ and $L^s$ 
if we do not tensor with ${\Bbb Z}[\frac{1}{2}]$. 

\begin{rem}{\rm The main ingredient behind the proof of 
Theorem \ref{mth} is [\cite{FJ}, theorem 
2.1 and remark 2.1.3]. This was also used before to prove 
the $\flfw$ and the $\flvw$ 
in \cite{FL} for elementary amenable groups and for torsion free 
virtually solvable subgroups of $GL(n, {\Bbb C})$.}\end{rem}

\begin{rem} {\rm We also note here that the $IC^{_{\ul L}{\cal H}^?_*}_{\cal {VC}}$ 
is known for many classes of groups. See [\cite{LR}, 5.3]. 
In \cite{HR} it was proved that the 
$IC^{_{L}{\cal H}^?_*}_{\cal {VC}}$ is true for the fundamental 
groups of closed manifolds with $\widetilde {\Bbb {SL}}\times 
{\Bbb E}^n$ structure for $n\geq 2$.}\end{rem}

\section{Statement of the Isomorphism Conjecture and some basic 
results in $L$-theory}
Given a normal subgroup $H$ of a group $G$ by [\cite{FR}, algebraic lemma] 
$G$ can be embedded as a subgroup in the wreath product $H\wr (G/H)$. We 
will always use this fact without mentioning. 

Let ${\cal H}^?_*$ be an equivariant homology theory with values in 
$R$-modules for $R$ a commutative associative ring with unit. In 
this article we are considering the special case $R={\Bbb Z}$.

In this section we always assume that a class of 
groups $\cal C$ is closed under 
isomorphisms, taking subgroups and taking quotients. We 
denote by ${\cal C}(G)$ the class of subgroups of a group 
$G$ which belong to $\cal C$. 
   
Given a group homomorphism $\phi:G\to H$ and $\cal C$ a class of 
subgroups of $H$ define $\phi^*{\cal C}$ by the class 
of subgroups $\{K<G\ |\ \phi (K)\in {\cal C}\}$ of $G$. For a class  
$\cal C$ of subgroups of a group $G$ there is a $G$-CW complex $E_{\cal 
C}(G)$ which is unique up to $G$-equivalence satisfying the property that 
for each $H\in {\cal C}$ the fixpoint set $E_{\cal C}(G)^H$ is 
contractible and $E_{\cal C}(G)^H=\emptyset$ for $H$ not in ${\cal C}$. 

\medskip
\noindent
{\bf (Fibered) Isomorphism Conjecture:} ([\cite{BL}, definition 1.1]) 
Let $G$ be a group 
and $\cal C$ be a class of subgroups of $G$. Then the {\it 
Isomorphism Conjecture} for the pair $(G, {\cal C})$ states that the 
projection 
$p:E_{\cal C}(G)\to pt$ to the point $pt$ induces an isomorphism $${\cal 
H}^G_n(p):{\cal H}^G_n(E_{\cal C}(G))\simeq {\cal H}^G_n(pt)$$ for $n\in 
{\Bbb Z}$. 

And the {\it Fibered Isomorphism Conjecture} for the pair $(G, {\cal 
C})$ states that for any group homomorphism $\phi: K\to G$ the 
Isomorphism Conjecture is true for the pair $(K, \phi^*{\cal C})$.

\begin{defn} \label{definition} ([\cite{R1}, definition 2.1]) 
{\rm Let $\cal C$ be a class of groups. If the (Fibered) 
Isomorphism Conjecture is
true for the pair $(G, {\cal C}(G))$ we say that the {\it (F)IC$^{{\cal
H}^?_*}_{\cal C}$ is true for $G$} or simply say {\it (F)IC$^{{\cal
H}^?_*}_{\cal C}(G)$ is satisfied}. Also we say that the {\it 
(F)ICwF$^{{\cal
H}^?_*}_{\cal C}(G)$ is satisfied} if
the {\it (F)IC$^{{\cal
H}^?_*}_{\cal C}$} is true for $G\wr H$ for any finite group
$H$.}\end{defn} 

Clearly, if $H\in {\cal C}$ then the (F)IC$^{{\cal
H}^?_*}_{\cal C}(H)$ is satisfied. 

The following is easy to prove and is known as 
the hereditary property of the Fibered Isomorphism 
Conjecture.

\begin{lemma} \label{heredi} If the FIC$^{{\cal H}^?_*}_{\cal C}(G)$ 
(the FICwF$^{{\cal H}^?_*}_{\cal 
C}(G)$) is satisfied then the 
FIC$^{{\cal H}^?_*}_{\cal C}(H)$ 
(the FICwF$^{{\cal H}^?_*}_{\cal C}(H)$) 
is satisfied for any subgroup $H$ of $G$.\end{lemma}

Let us denote by $_{\ul L}{\cal H}^?_*$ ($_{L}{\cal H}^?_*$), 
the equivariant homology 
theory arises for the $\ul L^{\l -\infty \r}$-theory 
($L^{\l -\infty \r}$-theory). In the statement of 
Theorem \ref{mth} the notation ${\cal H}^G_n(-, {\bf\ul L}^{\l -\infty \r})$ 
(${\cal H}^G_n(-, {\bf L}^{\l -\infty \r})$) 
stands for $_{\ul L}{\cal H}^G_n(-)$ ($_{L}{\cal H}^G_n(-)$),  
where ${\bf\ul L}^{\l -\infty \r}$ (${\bf L}^{\l -\infty \r}$)  
denotes the spectrum whose homotopy groups are the surgery 
groups ${\ul L}_*^{\l -\infty \r}$ (${L}_*^{\l -\infty \r}$). 
See [\cite{LR}, section 6.2] for details.

Throughout the article a `graph' is assumed to be connected 
and locally finite.

\begin{defn} \label{property} ([\cite{R1}, definition 2.2]) 
{\rm We say that ${\cal T}^{{\cal 
H}^?_*}_{\cal C}$ ($_w{\cal T}^{{\cal
H}^?_*}_{\cal C}$) is satisfied if for a
graph of groups $\cal G$ with vertex groups (and hence edge groups)
belonging to the class $\cal C$, 
the FIC$^{{\cal H}^?_*}_{\cal
C}$ (the FICwF$^{{\cal H}^?_*}_{\cal
C}$) for $\pi_1({\cal G})$ is true. 

We say that $_t{\cal T}^{{\cal H}^?_*}_{\cal C}$ 
($_{wt}{\cal T}^{{\cal H}^?_*}_{\cal C}$) 
is satisfied if for a graph of groups $\cal G$ with trivial edge
groups and the vertex
groups belonging to the class $\cal C$, the FIC$^{{\cal H}^?_*}_{\cal
C}$ (the FICwF$^{{\cal H}^?_*}_{\cal
C}$) for $\pi_1({\cal G})$ is true . 

And we say that {\it ${\cal P}^{{\cal H}^?_*}_{\cal C}$} is 
satisfied if for $G_1, G_2\in {\cal C}$ the product
$G_1\times G_2$ satisfies the FIC$^{{\cal H}^?_*}_{\cal C}$}\end{defn}

We start with the following general lemma which easily follows 
from [\cite{R1}, proposition 5.2] and 
[\cite{R2}, lemma 3.4].

\begin{lemma}\label{lemma3.4prop5.2} Assume that ${\cal P}^{{\cal 
H}^?_*}_{\cal 
C}$ is 
satisfied. 

(1). If the FIC$^{{\cal H}^?_*}_{\cal C}$ (the FICwF$^{{\cal 
H}^?_*}_{\cal C}$) is true 
for $G_1$ and $G_2$ then $G_1\times G_2$ 
satisfies the FIC$^{{\cal H}^?_*}_{\cal C}$ (the 
FICwF$^{{\cal
H}^?_*}_{\cal C}$). 

(2). Let $G$ be a finite index subgroup of a group $K$. If  
the group $G$ satisfies the FICwF$^{{\cal 
H}^?_*}_{\cal C}$ then $K$ also satisfies the 
FICwF$^{{\cal H}^?_*}_{\cal C}$.

(3). Let $p:G\to Q$ be a group homomorphism. If the FICwF$^{{\cal 
H}^?_*}_{\cal C}$ is true for $Q$ and for $p^{-1}(H)$ for all 
$H\in {\cal C}(Q)$ then the FICwF$^{{\cal
H}^?_*}_{\cal C}$ is true for $G$. 
\end{lemma} 
 
The following lemma is obvious.

\begin{lemma}\label{obvious} For a class of groups $\cal C$ 
the $\ficlu$ implies the 
$\ficl$ and the $\ficwlu$ implies 
the $\ficwl$.\end{lemma}
 
If the Isomorphism Conjecture is true for a group with respect to a class 
$\cal C$ of subgroups then obviously it is true for the group with 
respect to a class of subgroups containing $\cal C$. The following Lemma 
shows that sometimes the converse is also true.

\begin{lemma}\label{v-f} If a group $G$ satisfies the $\flvw$ 
(the $\flv$) then it also satisfies the $\flfw$ (the $\flf$).\end{lemma}

\begin{proof} See [\cite{FL}, lemma 5.1] or [\cite{LR}, proposition 2.18].\end{proof}

\begin{lemma} \label{limit} Let ${\cal C}={\cal {VC}}$ or $\cal {FIN}$. 
Let $\{G_i\}_{i\in I}$ be a directed sequence of
groups with
limit $G$. If each $G_i$ satisfies $\#$ 
then $G$ also satisfies $\#$. Here $\#$ is one of the 
followings.

1. $\ficl$. 2. $\ficlu$. 3. $\ficwl$. 4. $\ficwlu$.\end{lemma} 

\begin{proof} 
For $1$ and $2$ it directly follows from [\cite{FL}, theorem 7.1] and 
for $3$ and $4$ just note that for a finite group $F$, $G\wr F$ is a 
direct limit of $\{G_i\wr F\}_{i\in I}$ and then apply 
[\cite{FL}, theorem 7.1].\end{proof} 

\begin{lemma} \label{f-j} Let $\Gamma$ be a discrete 
cocompact subgroup of a Lie group
with finitely many connected components.
Then $\Gamma$ satisfies the $\flvuw$, $\flvw$, $\flfw$, $\flvu$, 
$\flv$ and the $\flf$.\end{lemma}

\begin{proof} By Lemma \ref{obvious} $\flvuw$ implies $\flvw$ and then 
apply Lemma \ref{v-f} to get $\flfw$. Similarly 
$\flvu$ implies $\flv$ and then applying Lemma \ref{v-f} we get 
$\flf$. Therefore we only have to show that $\Gamma$ satisfies 
$\flvuw$ and $\flvu$. For $\flvuw$ it follows directly 
from [\cite{FJ}, theorem 2.1 and remark 2.1.3]. 

For $\flvu$ we need the following lemma.

\begin{lemma} \label{lie} Let $G$ be a Lie group with finitely 
many components and 
let $F$ be a finite group. Then the wreath product $G\wr F$ is again 
a Lie group with finitely many components with respect to the product 
topology on $G^F\times F$ where $F$ is given the discrete 
topology and $G^F$ denotes the $|F|$-times direct product of 
$G$.\end{lemma}

\begin{proof} Recall that an element of $G^F$ is of the form 
$(g_{f_1},\ldots , g_{f_{|F|}})$ where $f_i\in F$. Now let $f\in F$. 
Then the regular action of $F$ on $G^F$ is by definition 
$f(g_{f_1},\ldots , g_{f_{|F|}})=(g_{f_1f^{-1}},\ldots , g_{f_{|F|}f^{-1}})$. 
It now follows from the definition of semi-direct 
product 
that the product and inverse operations on $G\wr F$ both are smooth. 
Therefore $G\wr F$ is a Lie group and clearly it has finitely many 
components.\end{proof} 

Now if $\Gamma$ is a discrete cocompact subgroup of $G$ and $G$ 
has finitely many components then it is easy to verify 
that $\Gamma\wr F$ is a discrete cocompact subgroup 
of $G\wr F$ for any finite group $F$. Here the Lie group structure 
on $G\wr F$ is as described in Lemma \ref{lie}. 

Hence we can again use [\cite{FJ}, theorem 2.1 and 
remark 2.1.3] to see that the $\flvu$ is satisfied for $\Gamma$.

This completes the proof.
\end{proof} 

\begin{lemma} \label{pvc} $\pvc$, $\pvcu$ and  $\pfu$ are 
satisfied.\end{lemma}

\begin{proof} Recall that $\pvc$ states that the 
$\flvuw$ is true for $V_1\times V_2$ for any two virtually 
cyclic groups $V_1$ and $V_2$. Let $V_1$ and $V_2$ be two 
such groups then $V_1\times V_2$ contains a free abelian 
normal subgroup $H$ (say) (on at most $2$ generators) of 
finite index. Hence $V_1\times V_2$ is a subgroup of 
$H\wr ((V_1\times V_2)/H)$. Therefore by Lemma \ref{heredi} it 
is enough to prove the $\flvuw$ for $H\wr ((V_1\times V_2)/H)$.  

Now we need the following well known fact.

\begin{lemma}\label{surface} Let $S$ be a closed orientable surface of 
genus $\geq 1$. Then $\pi_1(S)$ is a discrete 
cocompact subgroup of a Lie group with finitely many 
components.\end{lemma}

\begin{proof} If the genus of $S$ is $1$ then $\pi_1(S)$ is a 
discrete cocompact subgroup of the Lie group of isometries of the 
flat Euclidean plane. And 
if the genus of $S$ is $\geq  2$ then the corresponding 
Lie group is the group of isometries of the hyperbolic 
plane.\end{proof} 

If $V_1\times V_2$ is virtually cyclic then there is 
nothing to prove. If $H$ has rank $2$ then 
applying Lemmas \ref{surface} and 
\ref{f-j} we see that $\pvc$ is satisfied. 
Next we apply Lemma \ref{obvious} to see that 
$\pvcu$ is also satisfied. And $\pfu$ is obvious.\end{proof}

\begin{defn} {\rm Choose two classes of groups ${\cal C}_1$ and 
${\cal C}_2$ so that ${\cal C}_1\subset {\cal C}_2$. 
We say that ${\cal T}^{{\cal H}^?_*}_{{\cal C}_1, {\cal C}_2}$ 
is satisfied if for a graph of groups $\cal G$ 
with vertex groups belonging to the class ${\cal C}_1$ the 
$FIC^{{\cal H}^?_*}_{{\cal C}_2(\pi_1({\cal G}))}
(\pi_1({\cal G}))$ is satisfied.}\end{defn}

\begin{lemma}\label{tvc} $\tfinv$, 
$\tfinvu$, $\wtfinu$, $\wttv$ and $\wttvu$ are
satisfied.\end{lemma}

\begin{proof} At first we check $\tfinv$. So let $G$ be a group 
and $\cal G$ be a  
graph of finite groups with $\pi_1({\cal G})\simeq G$. If $\cal G$ 
is an infinite graph then we write $\cal G$ as an increasing union 
of finite subgraphs ${\cal G}_i$. Then $\pi_1({\cal G})\simeq 
\lim_{i\to\infty}\pi_1({\cal G}_i)$. Hence using Lemma \ref{limit} 
we can assume that $\cal G$ is finite. It is now well known 
that $\pi_1({\cal G})$ contains a finitely generated free 
subgroup of finite index. See [\cite{R1}, lemma 3.2].   
$\tfinv$ now follows from the following Main Lemma.

To check $\tfinvu$ and $\wtfinu$ we only need to use 
Lemmas \ref{obvious} and \ref{v-f}. 

Next we prove $\wttv$ and 
$\wttvu$. 

Let $\cal G$ be a graph of groups with virtually 
cyclic vertex groups and trivial edge groups. As before we can assume that 
$\cal G$ is finite. Hence the group $\pi_1({\cal G})$ is 
virtually free. This follows from the 
following two Lemmas. 

\begin{lemma}\label{action} If a graph of groups $\cal G$ has 
trivial edge groups then $\pi_1({\cal G})$ is isomorphic to 
the free product of a free group and the vertex 
groups of $\cal G$.\end{lemma}

\begin{proof} Apply [\cite{R1}, lemma 6.2] and note that 
$\cal G$ is a direct limit of its finite subgraphs 
of groups.\end{proof}

\begin{lemma} Let $V_1$ and $V_2$ be two virtually free  
groups then $V_1*V_2$ is virtually free.\end{lemma}

\begin{proof} We have a surjective homomorphism 
$p:V_1*V_2\to V_1\times V_2$. Let $H_i$ be a free 
subgroup of $V_i$ of finite index for $i=1,2$. Hence 
$H=H_1\times H_2$ has finite index in $V_1\times V_2$. Note that 
$V_1*V_2$ acts on a tree with trivial edge stabilizers 
and the vertex stabilizers are conjugate to $V_1$ or $V_2$. 
Hence $p^{-1}(H)$ also acts on the same tree. It  
follows that the edge stabilizers of this restricted 
action are again trivial and the vertex stabilizers 
are free. And hence $p^{-1}(H)$ is a free group 
by Lemma \ref{action}. This completes the 
proof.\end{proof}

Therefore we can apply Lemmas 
\ref{limit} and the following Main Lemma to complete the proof of 
Lemma \ref{tvc}.
\end{proof}  

\begin{lemma}\label{dproduct} Assume that the $\#$ is true for two 
groups $G_1$ and $G_2$ then the $\#$ is true 
for the direct product $G_1\times G_2$. Here $\#$ denotes one of the 
followings.

1. $\flvuw$. 2. $\flvw$. 3. $\flvu$. 4. $\flv$.

\end{lemma} 

\begin{proof} The proof is a combination of Lemma \ref{pvc} 
and 
$(1)$ of Lemma \ref{lemma3.4prop5.2}.\end{proof}

\begin{mainl} The $\flvu$, $\flv$ and $\flf$ are true for any 
virtually free group.\end{mainl}

\begin{proof} We only prove the lemma for the $\flvu$. The other 
two conclusions will follow using Lemmas \ref{obvious} and \ref{v-f}.

Let $\Gamma$ be a virtually free group and $G$ be a
free normal subgroup of $\Gamma$ with $F$ the finite quotient group. Let
$F'$ be another finite group and denote by $C$ the wreath product $F\wr
F'$. Then we have the following inclusions.
$$\Gamma\wr F' < (G\wr F)\wr F' < G^{F\times F'}\wr C < (G\wr
C)\times\cdots \times (G\wr C).$$ 

The inclusions are easy to check. (See [\cite{R0}, lemma 5.4] for 
the second inclusion). There are
$|F\times F'|$ factors in the last term. Therefore using Lemmas 
\ref{heredi} and \ref{dproduct} we see that it 
is enough to
prove the $\flvuw$ for $G\wr C$ for an arbitrary finite group $C$. 
Equivalently, we need to prove the $\flvu$ for $G$.
If $G$
is infinitely generated then let $G$ be the limit of a 
sequence of finitely
generated subgroups of $G$. 
Hence by Lemma \ref{limit} we can assume that $G$ is
finitely generated. Therefore $G$ is isomorphic to the 
fundamental group of an orientable $2$-manifold $M$ with boundary. 
Consequently, 
$G$ is isomorphic to a subgroup of $\pi_1(M\cup_{\partial} M)$,
where $M\cup_{\partial} M=S$ (say) denotes the double of $M$. Again using
Lemma \ref{heredi} it is enough to prove the $\flvu$ for $\pi_1(S)$, 
where $S$ is a closed orientable surface. Without 
loss of generality we
can assume that $S$ has genus $\geq 1$. 
Now applying Lemmas \ref{surface} and \ref{f-j} 
we complete the proof of the 
Main Lemma.\end{proof} 

\begin{prop} \label{fproduct} Assume that the $\#$ is true for two
groups $G_1$ and $G_2$ then the $\#$ is true
for the free product $G_1*G_2$ also. Here $\#$ is 
one of the followings: 

1. $\flvuw$. 2. $\flvw$. 3. $\flfw$. 4. $\flvu$. 
5. $\flv$. 6. $\flf$.
\end{prop}

\begin{proof} The proof follows from Lemmas \ref{pvc}, \ref{tvc} and 
[\cite{R1}, lemma 6.3].\end{proof}

\begin{lemma} \label{import} Let $1\to K\to G\to N\to 1$ 
be an exact sequence of groups. 
Then the followings hold.

1. If the $\ficwf(K)$ and the $\ficf(N)$ are satisfied 
then the $\ficf(G)$ is also satisfied.

2. If the $\ficwf(K)$ and the $\ficwf(N)$ are satisfied 
then the $\ficwf(G)$ is also satisfied.
\end{lemma} 

\begin{proof} Apply $(2)$ and $(3)$ of Lemma \ref{lemma3.4prop5.2} and 
note that ${\cal P}^{{\cal H}^?_*}_{\cal {FIN}}$ is satisfied.\end{proof}

\section{Braid groups}
Let ${\Bbb C}^N$ be the $N$-dimensional complex space. A {\it hyperplane 
arrangement} in ${\Bbb C}^N$ is by definition a finite collection 
$\{V_1, V_2,\ldots , V_n\}$ of $(N-1)$-dimensional linear subspaces of
${\Bbb C}^N$. 

Now we recall the definition of a fiber-type 
hyperplane arrangement from [\cite{OT}, page 162]. Let us denote by
${\cal V}_n$ the arrangement $\{V_1, V_2,\ldots , V_n\}$ in ${\cal
C}^N$. ${\cal V}_n$ is called {\it strictly linearly fibered} 
if after a suitable linear change of coordinates, the restriction of the 
projection of ${\Bbb C}^N-\cup_{i=1}^n V_i$ to the first $(N-1)$
coordinates 
is a fiber bundle projection whose base space is the complement of an 
arrangement ${\cal W}_{n-1}$ in ${\Bbb C}^{N-1}$ and whose fiber is the 
complex plane minus finitely many points. By definition the arrangement 
$0$ in ${\Bbb C}$ is fiber-type and ${\cal V}_n$ is defined to be 
{\it fiber-type} if ${\cal V}_n$ is strictly linearly fibered 
and ${\cal W}_{n-1}$ is of fiber type. It follows by repeated 
application 
of homotopy exact sequence for fibration that the complement 
${\Bbb C}^N-\cup_{i=1}^n V_i$ is aspherical and hence the fundamental 
group is torsion free.

\begin{lemma} ([\cite{FR}, theorem 5.3]) 
$\pi_1({\Bbb C}^N-\cup_{i=1}^n V_i)$ 
is a strongly poly-free group.\end{lemma}

Now recall that the pure braid group $PB_n$ on $n$ strings 
is by definition $\pi_1({\Bbb C}^{n+1}-\cup_{i,j} V_{ij})$ where 
$V_{ij}$ is the hyperplane $x_i=x_j$ for $i<j$ and $x_i$'s being  
the coordinates in ${\Bbb C}^{n+1}$. One can show that 
$\{V_{ij}\}$ is a fiber-type arrangement and hence 
$PB_n$ is a strongly poly-free group. 
See [\cite{AFR}, theorem 2.1]. 

The full braid group
$B_n$ is by definition 
$\pi_1(({\Bbb C}^{n+1}-\cup_{i,j} V_{ij})/S_{n+1})$ where the 
symmetric group $S_{n+1}$ on $(n+1)$-symbols acts on 
${\Bbb C}^{n+1}-\cup_{i,j} V_{ij}$ by permuting the coordinates. 
This action is free and therefore $PB_n$ is a normal 
subgroup of $B_n$ with quotient $S_{n+1}$. 
 
Recall that in \cite{AFR} we proved the following.

\begin{thm} (theorem 1.3 and corollary 1.4 in \cite{AFR})
\label{afr} Let
$\Gamma$ be the fundamental group of a fiber-type hyperplane arrangement
complement or more generally a strongly poly-free group. Then
$Wh(\Gamma)=\tilde K_0({\Bbb Z}\Gamma)=K_i({\Bbb Z}\Gamma)=0$ for
$i<0$.\end{thm}

Theorem \ref{afr} and an application of the Rothenberg's exact sequence
show the following. See [\cite{R}, 17.2].

\begin{lemma} \label{rothen} Let $\Gamma$ be as in Theorem
\ref{afr} then $L^{\l -\infty \r}_n({\Bbb Z}\Gamma)\simeq L^h_n({\Bbb
Z}\Gamma)\simeq L^s_n({\Bbb Z}\Gamma)$.\end{lemma} 

In the situation of
$\Gamma$ as in Theorem \ref{afr}, 
Lemma \ref{rothen} shows that the $2$-torsions which 
appear in the three surgery
groups are also isomorphic.

In \cite{R3} we will use Corollary \ref{classic}, 
Theorem \ref{hyper} and Lemma \ref{rothen} to calculate 
explicitly the surgery groups 
of the fundamental groups of fiber-type hyperplane arrangement 
complement in the complex space. In particular this applies to 
the pure braid groups. 

\section{Proof of Theorem \ref{hyper}}
\begin{proof}[Proof of Theorem \ref{hyper}]
{\bf 1 \& i.} Since $\uF$ ($\F$) contains the 
discrete cocompact subgroups 
of (linear) Lie groups with finitely many components it follows that 
finite groups and the infinite cyclic group 
belong to $\uF$ ($\F$). Next apply the `polycyclic extension' 
(`extension') condition to 
complete the proof of $(1)$ ($(i)$).

\medskip

\noindent
{\bf 2 \& ii.} At first note that a countable 
infinitely generated  
group is a direct limit of finitely generated subgroups.

Now using $(1)$ ($(i)$) and 
the `free product' condition we get that 
finitely generated free groups belong to $\F$ ($\uF$) and since 
$\F$ ($\uF$) has the property `direct limit' the proof 
follows for infinitely generated free groups. 

Using $(1)$ ($(i)$) and the `extension' (`polycyclic extension') 
condition we see that 
finitely generated abelian groups belong to $\F$ ($\uF$). 
Therefore countable abelian groups belong 
to $\F$ ($\uF$) by the `direct limit' condition.

\medskip

\noindent
{\bf 3.} The 
proof is by induction on the index of the poly-free group. 
If $n=1$ then $G$ is free and 
hence $G\in \F$ (apply $(2)$). So assume that poly-free 
groups of index $\leq n-1$ belong to $\F$ and let 
$G$ has an index $n$ filtration. Note that $G_{n-1}$ 
is a poly-free group of index $n-1$ and $G/G_{n-1}$ 
is a free group. Since $\F$ is closed under extensions 
we can now apply $(2)$ and 
the induction hypothesis to show that $G\in\F$.

\medskip

\noindent
{\bf 4 \& 5.} Recall from Section 3 that pure braid groups 
are strongly poly-free and strongly poly-free groups 
are poly-free. Also the full braid group $B_n$ contains 
the pure braid group $PB_n$ as a subgroup of finite index. 
Hence using {\bf 2} and since $\F$ is closed under extensions 
the proofs of $3$ and $4$ are complete.

\medskip

\noindent
{\bf 6.} Let $\G$ be a cocompact discrete subgroup 
of a Lie group with finitely many components. Then following 
the steps in the proof of [\cite{R2}, $2(a)$ of theorem 2.2] or 
of [\cite{FJ}, theorem 2.1]  
we have the following three exact sequences.

$$1\to F \to \G \to \G' \to 1.$$

$$1\to \G_R \to \G' \to \G_S \to 1.$$

$$1\to \G_H \to \G_S \to \G_L \to 1.$$

Where $F$ is finite, $\G_R$ is virtually poly-${\Bbb Z}$, $\G_H$ 
is virtually finitely generated abelian and $\G_L$ is 
a cocompact discrete subgroup of a Linear Lie group with 
finitely many components. Now note that finitely generated 
free abelian groups and poly-${\Bbb Z}$ 
groups (since they are also poly-free) belong 
to $\F$. Therefore, we can again apply the hypothesis 
that $\F$ is closed under extensions and use the above three exact 
sequences to complete the proof of ${\bf 5}$.

\medskip

\noindent
{\bf 7.} Let $\G$ be a group so that $\G^{(n)}\in\F$ for 
some $n$. Using the extension condition it is enough to 
show that $\G/\G^{(n)}\in\F$, that is we need to show 
that $\F$ contains the solvable groups.  
 
So let $\G$ be a solvable group. 
Using the `direct limit' condition in the definition 
of $\F$ we can assume that 
$\G$ is finitely generated, for any countable infinitely 
generated group is a direct limit of finitely generated 
subgroups.

We say that $\G$ is $n$-step 
solvable if $\G^{(n+1)}=(1)$ and $\G^{(n)}\neq (1)$. The 
proof is by induction on $n$. Since 
countable abelian groups belong to $\F$ 
(by $(2)$), the 
induction starts. 

So assume that a finitely generated $k$-step solvable 
group for $k\leq n-1$ belongs to $\F$ and $\G$ is 
$n$-step solvable.

We have the following exact

sequence. $$1\to \G^{(n)}\to \G\to \G/\G^{(n)}\to 1.$$

Note that $\G^{(n)}$ is abelian and $\G/\G^{(n)}$ is 
$(n-1)$-step solvable. Using the `extension' condition 
and the induction hypothesis we complete the proof.

\medskip

\noindent
{\bf iii.} This follows from the definition of $\uF$.

\medskip

\noindent
{\bf iv.} Using the `polycyclic extension' condition and the 
following Lemma we complete the proof.

\begin{lemma} \label{wikling} Let $G$ be a 
virtually polycyclic group. Then 
$G$ contains a finite normal subgroup so that the 
quotient is a discrete cocompact subgroup of a 
Lie group with finitely many components.\end{lemma}

\begin{proof} See [\cite{W}, theorem 3, remark 4 
in page 200].\end{proof} 

\end{proof}

\section{Proofs of Theorems \ref{mth} and \ref{mth1}} 
\begin{proof}[Proof of Theorem \ref{mth}]
It is enough to prove 
the followings. 

\noindent
{\bf A.} The $\flf$ ($\flvu$) for the groups appearing in $(1)$ 
($(i)$) of  
Definition \ref{begin} is true. 

\noindent
{\bf B.} The 
statement {\it `The $\flf$ ($\flvu$) is satisfied'} 
is closed under the operations 
described in $2$ to $5$ ($(ii)$ to $(iv)$) of 
Definition \ref{begin}. 

\medskip

\noindent
{\bf Proof of A.} 
The proof follows from Lemma \ref{f-j}.

\medskip

\noindent
{\bf Proof of B.} The proof follows from Theorem \ref{mth1}.

This proves the Theorem \ref{mth}.\end{proof}

\begin{proof}[Proof of Theorem \ref{mth1}] For $\cal C$: 
2, 3, 4 and 5 follows from Lemma \ref{heredi}, Proposition \ref{fproduct}, 
Lemma \ref{limit} and Lemma \ref{import} respectively.

For $\cal D$: $(ii)$ follows from Lemma \ref{heredi}, Proposition 
\ref{fproduct} and Lemma \ref{limit}. $(iii)$ follows 
from Lemma \ref{dproduct}. $(iv)$ follows using $(3)$ of 
Lemma \ref{lemma3.4prop5.2}, Lemma \ref{wikling} and Lemma 
\ref{f-j}. To apply $(3)$ of Lemma \ref{lemma3.4prop5.2} we
will need the fact that if a group contains a finite normal 
subgroup with virtually cyclic quotient then the group is 
virtually cyclic. This follows from [\cite{R1}, lemma 6.1].
\end{proof}

\newpage
\bibliographystyle{plain}
\ifx\undefined\bysame
\newcommand{\bysame}{\leavevmode\hbox to3em{\hrulefill}\,}
\fi

\medskip

\end{document}